\documentclass{amsproc}
\usepackage{amssymb}
\usepackage{graphicx}
\usepackage{amscd}
\usepackage{amsmath}
\usepackage{showkeys}
\theoremstyle{plain}

\newtheorem{definition}{Definition}

\newtheorem{lemma}{Lemma}

\numberwithin{equation}{section}

\newcommand{\N}{\mathbb{N}}
\newcommand{\p}{\partial}

\begin{document}

\title{Spectral Estimation Problem in Infinite Dimensional Spaces}
\author{ S. A. Avdonin}
\address{Department of Mathematics
and Statistics University of Alaska Fairbanks, PO Box 756660
Fairbanks, AK 99775, USA.} \email{s.avdonin@alaska.edu}

\author{ V. S. Mikhaylov}
\address{St. Petersburg   Department   of  the V.A. Steklov    Institute   of   Mathematics
of   the   Russian   Academy   of   Sciences, 7, Fontanka, 191023
St. Petersburg, Russia} \email{vsmikhaylov@pdmi.ras.ru}
\thanks{The research was supported by RFFI 11-01-00407A, RFFI 12-01-31446}
\date{November 30, 2012}

\maketitle

\noindent {\bf Abstract.} We consider the generalized spectral
estimation problem in infinite dimensional spaces. We solve this problem using the boundary
control approach to inverse theory and provide an application to the initial boundary value problem for
a hyperbolic system.

\section{Introduction}

The classical spectral estimation problem consists of the recovery
of the coefficients $a_n$, $\lambda_k$, $k = 1,\ldots,N, \; N \in
\mathbb{N},$ of a signal
\begin{equation*}
s(t)=\sum_{n=1}^N a_ke^{\lambda_k t},\quad t\geqslant 0
\end{equation*}
from the given observations $s(j)$, $j=0,\ldots, 2N-1,$ where the
coefficients $a_k,$ $\lambda_k$ may be arbitrary complex numbers.
The literature describing variuos methods for solving the spectral
estimation problem is very extensive: see for example the list of
references in \cite{AGM}.

In papers \cite{AB,ABN}  a new approach to
this problem was proposed. In this approach the signal $s(t)$ was treated as a
kernel of certain convolution operator corresponding to an
input-output map for some linear discrete-time dynamical system.
While the system realized from the input-output map is not unique,
the coefficients $a_n$ and $\lambda_n$ can be
determined uniquely using the non-selfadjoint version of the boundary
control method \cite{ABe}.

Later on the infinite-dimensional version of this method has been
developed in \cite{AGM}. More precisely, the problem of the
recovering the coefficients $a_k, \lambda_k\in \mathbb{C}$, $k\in
\mathbf{N},$ of the given signal
\begin{equation*}
s(t)=\sum_{k=1}^\infty a_ke^{\lambda_kt},\quad t\in (0,T),
\end{equation*}
provided the sum converges in $L_2(0,T)$ was solved there.  In the
present paper we solve the so-called generalized   spectral
estimation problem. It is set up in the following way: to recover
the coefficients $a_k(t)$, $\lambda_k$, $k\in \mathbb{N},$ of a
signal
\begin{equation}
\label{unknown_f} S(t)=\sum_{k=1}^\infty a_k(t)e^{\lambda_k
t},\quad t\in (0,T),
\end{equation}
from the given data $S\in L_2(0,T).$ We assume that $T$ ia a
positive number,  $\lambda_k\in \mathbb{C}$ and for each $k,$
$a_k(t)=\sum_{i=0}^{L_k}a^i_kt^i$ are polynomials of the order
$L_k$ with complex valued coefficients $a^i_k.$

In Section 2 we  recover the unknown parameters $\, \lambda_k, \,
L_k, \, a^i_k \,; $ $\, i=0,\ldots,L_k, \; k \in \mathbb{N},$ from
$S(t)$, $t\in (0,T)$. In Section 3, as an application of the
generalized spectral estimation, we consider the continuation
problem of the inverse dynamical data in the identification
problem for the first order hyperbolic system.

\section{Spectral estimation. The case of multiple poles}

We consider the dynamical system in a complex Hilbert space $H$:
\begin{equation}
\label{dyn_sys} \dot x(t)=Ax(t)+bf(t), \quad t\in (0,T),\quad
x(0)=0.
\end{equation}
Here $b\in H$, $f\in L_2(0,T),$ and we assume that the spectrum of
the operator $A$, $\{\lambda_k\}_{k=1}^\infty$ is not simple. We denote the
algebraic multiplicity of $\lambda_k$  by $L_k,$
$k \in \mathbb{N},$ and assume also that the set of all root
vectors $\{\phi_k^i\},$
$i=1,\ldots, L_k,$  $k\in \mathbb{N},$ forms a Riesz basis in $H$. Here the
vectors from the chain $\{\phi^i_k\}_{i=1}^{L_k}$,
$k \in \mathbb{N},$ satisfy the equations
\begin{equation*}
\left(A-\lambda_k\right)\phi_k^1=0,\quad
\left(A-\lambda_k\right)\phi_k^i=\phi_k^{i-1},\,\, 2\leqslant
i\leqslant L_k.
\end{equation*}
Along with (\ref{dyn_sys}), we consider the dynamical system for
the adjoint operator:
\begin{equation}
\label{dyn_sys_adj} \dot y(t)=A^*y(t)+dg(t), \quad t\in
(0,T),\quad y(0)=0,
\end{equation}
where $d\in H$, $g\in L_2(0,T)$.  The spectrum of $A^*$ is
$\{\overline\lambda_k\}_{k=1}^\infty$ and the root vectors
$\{\psi_{k}^{i}\}_{i=1}^{L_k}$, $i=1,\ldots, L_k,$  $k\in \mathbb{N},$  also form a
Riesz basis in $H$ and satisfy the equations
\begin{equation*}
\left(A^*-\overline\lambda_k\right)\psi_{k}^{L_k}=0,\quad
\left(A^*-\overline\lambda_k\right)\psi_{k}^i=\phi_{k}^{i+1},\,\,
1\leqslant i\leqslant L_k-1.
\end{equation*}
Moreover, the root vectors of $A$ and $A^*$ are normalized in
accordance with
\begin{eqnarray*}
\left<\phi_k^i,\psi_{l}^j\right>=0, \text{ if $k\not=l$ or
$i\not= j$};\\
\left<\phi_k^i,\psi_{k}^i\right>=1,\,\,\,i=1,\ldots,L_k, \, k\in
\mathbb{N}.
\end{eqnarray*}
We  consider $f$ and $g$ as the inputs of  the systems (\ref{dyn_sys}) and (\ref{dyn_sys_adj}) and define the outputs $z$ and $w$ by the formulas
\begin{equation*}
%\label{responce}
z(t)=\left<x(t),d\right> ,\quad
w(t)=\left<y(t),b\right>.
\end{equation*}
Suppose that the vector $b$ has a representation
$b=\sum_{k=1}^\infty\sum_{i=1}^{L_k}b_k^i\phi_{k}^i$. We look for
the solution to (\ref{dyn_sys}) in the form
\begin{equation}
\label{dyn_sys_sol_repr}
x(t)=\sum_{k=1}^\infty\sum_{i=1}^{L_k}c_k^i(t)\phi_{k}^i.
\end{equation}
Plugging (\ref{dyn_sys_sol_repr}) into (\ref{dyn_sys}),
multiplying  by $\psi_{k}^i$,
$i=1,\ldots,L_k$,  $k \in \mathbb{N}$, we get the following equations for $c_k^i(t)$:
\begin{eqnarray*}
&\dot c_k^{L_k}(t)=\lambda_k c_k^{L_k}(t)+b_k^{L_k}f(t),\quad c_k^{L_k}(0)=0,\\
&\dot c_k^{i}(t)=\lambda_k
c_k^{i}(t)+c_k^{i+1}(t)+b_k^{i}f(t),\,\, c_k^{i}(0)=0,\,\,
i=1,\ldots,L_k-1.
\end{eqnarray*}
%Multiplying by $\phi_{k*}$, $k=N+1,\ldots\infty$, we have
%\begin{equation}
%\label{od3}
%\dot c_k(t)=\lambda_k c_k(t)+b_kf(t),\quad c_k(0)=0.
%\end{equation}

Solving the system of ODEs  we find the coefficients  $c_k^i(t)$:
\begin{eqnarray*}
c_k^{L_k}(t)=\int_0^te^{\lambda_k(t-\tau)}b_k^{L_k}f(\tau)\,d\tau,\\
c_k^{L_k-1}(t)=\int_0^te^{\lambda_k(t-\tau)}\left[(t-\tau)b_k^{L_k}+b_k^{L_k-1}\right]f(\tau)\,d\tau,\\
c_k^{L_k-2}(t)=\int_0^te^{\lambda_k(t-\tau)}\left[\frac{(t-\tau)^2}{2}b_k^{L_k}+(t-\tau)b_k^{L_k-1}+b_k^{L_k-2}\right]f(\tau)\,d\tau,\\
c_k^{1}(t)=\int_0^te^{\lambda_k(t-\tau)}\left[\frac{(t-\tau)^{L_k-1}}{(L_k-1)!}b_k^{L_k}+\ldots+(t-\tau)b_k^{2}+b_k^{1}\right]f(\tau)\,d\tau.
\end{eqnarray*}
%and solving (\ref{od3}) we obtain for $k=N+1,\ldots\infty$
%\begin{equation*}
%c_k(t)=\int_0^te^{\lambda_k(t-\tau)}b_kf(\tau)\,d\tau.
%\end{equation*}

Similarly, we represent the vector $d$ in the form
$d=\sum_{k=1}^\infty\sum_{i=1}^{L_k}d_k^i\psi_{k}^i.$ Then the output $z$ can be written as
\begin{equation*}
%\label{responce_int}
z(t)=\left<x(t),d\right>=\sum_{k=1}^\infty\sum_{i=1}^{L_k}c_k^i(t)d_k^i=\int_0^tr(t-\tau)f(\tau)\,d\tau,
\end{equation*}
where $r(t)$ is defined as
\begin{equation}
\label{r_representation} r(t)=\sum_{k=1}^\infty e^{\lambda_k
t}\left[a_k^1+a_k^2t+a_k^3\frac{t^2}{2}+\ldots
+a_k^{L_k-1}\frac{t^{L_k-2}}{(L_k-2)!}+a_k^{L_k}\frac{t^{L_k-1}}{(L_k-1)!}\right].
\end{equation}
Here we introduced the notations
\begin{eqnarray}
a_k^1=\sum_{i=1}^{L_k}b_k^i d_k^i, \quad
a_k^2=\sum_{i=2}^{L_k}b_k^i d_k^{i-1},\quad
a_k^3=\sum_{i=3}^{L_k}b_k^i d_k^{i-2},\ldots\label{coeff_repr}\\
\ldots a_k^{L_k-1}=\sum_{i=L_k-1}^{L_k}b_k^i d_k^{i-(L_k-2)},\quad
a_k^{L_k}=b_k^{L_k} d_k^{1},\quad k \in \mathbb{N}.\notag
%\\
%a_k=b_kd_k,\quad k=1,\ldots,\infty.\label{coeff_repr_3}
\end{eqnarray}

It is important to notice notice that the {\sl response function}  $r(t)$ has the form of the series in
(\ref{unknown_f}).

Looking for the solution of (\ref{dyn_sys_adj}) in the form
\begin{equation*}
y(t)=\sum_{k=1}^\infty\sum_{i=1}^{L_k}h_k^i(t)\psi_{k}^i,
\end{equation*}
we derive the following system of ODEs for $h_k^i(t)$, $i=1,\ldots,L_k$, $k \in \mathbb{N}$:
\begin{eqnarray*}
&\dot h_k^{1}(t)=\overline\lambda_k h_k^{1}(t)+d_k^{1}g(t),\quad h_k^{1}(0)=0,\\
&\dot h_k^{i}(t)=\overline\lambda_k
h_k^{i}(t)+h_k^{i-1}(t)+d_k^{i}g(t),\,\, h_k^{i}(0)=0,\,\,
i=2,\ldots,L_k.
\end{eqnarray*}
%For  $k=1,\ldots\infty$ we derive
%\begin{equation}
%\dot h_k(t)=\overline\lambda_k h_k(t)+d_kg(t),\quad
%h_k(0)=0.\label{od_ad_3}
%\end{equation}
Solving this system we obtain the coefficients $h_k^i(t)$:
\begin{eqnarray*}
h_k^{1}(t)=\int_0^te^{\overline\lambda_k(t-\tau)}d_k^{1}g(\tau)\,d\tau,\\
h_k^{2}(t)=\int_0^te^{\overline\lambda_k(t-\tau)}\left[(t-\tau)d_k^{1}+d_k^{2}\right]g(\tau)\,d\tau,\\
h_k^{3}(t)=\int_0^te^{\overline\lambda_k(t-\tau)}\left[\frac{(t-\tau)^2}{2}d_k^{1}+(t-\tau)d_k^{2}+d_k^{3}\right]g(\tau)\,d\tau,\\
h_k^{L_k}(t)=\int_0^te^{\overline\lambda_k(t-\tau)}\left[\frac{(t-\tau)^{L_k-1}}{(L_k-1)!}d_k^{1}+\ldots+(t-\tau)d_k^{L_k-1}+d_k^{L_k}\right]g(\tau)\,d\tau,
\end{eqnarray*}
%Solving (\ref{od_ad_3}) for  $k=N+1,\ldots\infty$ yields
%\begin{equation*}
%h_k(t)=\int_0^te^{\overline\lambda_k(t-\tau)}d_kg(\tau)\,d\tau,
%\end{equation*}
The output of the system (\ref{dyn_sys_adj}) is given by
\begin{equation*}
w(t)=\left<y(t),b\right>=\sum_{k=1}^\infty\sum_{i=1}^{L_k}h_k^i(t)b_k^i=
\int_0^t \overline{r(t-\tau)}g(\tau)\,d\tau.
\end{equation*}

We introduce now the {\sl connecting operator} $C^T: L_2(0,T)\mapsto
L_2(0,T)$ defined through its bilinear form by the formula:
\begin{equation*}
\left<C^Tf,g\right>=\left<x(T),y(T)\right>.
\end{equation*}
\begin{lemma}
The connecting operator $C^T$ has a representation
$(C^Tf)(t)=(Rf)(2T-t)$, or
\begin{equation*}
(C^Tf)(t)=\int_0^T r(2T-t-\tau)f(\tau)\,d\tau.
\end{equation*}
\end{lemma}
\begin{proof}
We introduce the function
$\chi(s,t):=\left(x(s),y(t)\right)_H$. It is straightforward to
check that for $s,t>0$, this function satisfies the equation
\begin{equation*}
\chi_t(s,t)-\chi_s(s,t)=(r*f)(s)g(t)-(r*g)(t)f(s)
\end{equation*}
with the boundary conditions $\chi(0,t)=\chi(s,0)=0$.
This initial boundary value problem can be solved explicitly.
Since $x(T)$ and $y(T)$ are independent of the value of $f(t)$ and $g(t)$ for $t>T,$
we may put $f(t)=g(t)=0,$ if  $t>T,$ when compute $(C^Tf,g)_H.$
Taking this into account, we obtain:
\begin{eqnarray*}
\left<C^Tf,g\right>=\chi(T,T)=
\int_0^T\int_0^{2T-\gamma}r(2T-\gamma-\tau)f(\tau)
g(\gamma)\,d\tau \,d\gamma,
\end{eqnarray*}
and therefore,
\begin{equation} \label{ct}
(C^Tf)(t)=\int_0^{2T-t} r(2T-t-\tau)f(\tau)\,d\tau=\int_0^{T}
r(2T-t-\tau)f(\tau)\,d\tau.
\end{equation}
\end{proof}

Next, we demonstrate how to find $\lambda_k$, $L_k$ and
$a_k^i$,  $i=1,\ldots,L_k, $ $k \in \mathbb{N},$  given the function $r(t)$
in the form (\ref{unknown_f}).
To do that we use
 the ideas of the boundary control method, more precisely,  the
possibility to extract the spectral data from the dynamical data
(see \cite{B07,B03}). We assume that the system (\ref{dyn_sys}) is
spectrally controllable in time $T$. This means that, for any $i \in \{1,\ldots,L_k\}$
and  any $k \in \mathbb{N},$  there
exists  $\{f_k^i\}\in H^1_0(0,T)$,
 such that $x^{f_k^i}(T)=\phi_k^i$. By the definition of $\{f_k^i\},$
\begin{equation}
\label{i_1} \dot
x^{f_k^1}(T)=Ax^{f_k^1}(T)+bf_k^1(T)=A\phi_k^1=\lambda_k\phi_k^1=\lambda_k
x^{f_k^1}(T),\,\, k \in \N,
\end{equation}
\begin{equation}
\label{i_2} \dot
x^{f_k^i}(T)=A\phi_k^i=\lambda_k\phi_k^i+\phi_k^{i-1}=\lambda_k
x^{f_k^i}(T)+x^{f_k^{i-1}}(T),\,\, i=2,\ldots,L_k, \,\, k \in \N.
%\\
%x^{f_k}(T)=Ax^{f_k}(T)+bf_k(T)=A\phi_k=\lambda_k\phi_k=\lambda_k
%x^{f_k}(T),\,\, k=N+1,\ldots\infty.\label{i_3}
\end{equation}
The definition of the operator $C^T$ and equations (\ref{i_1}) imply that for any
 $g \in L_2(0,T)$ one has
\begin{eqnarray*}
\left<C^T\dot f_k^1,g\right>=\left<x^{\dot
f_k^1}(T),y^g(T)\right>=\left<\dot
x^{f_k^1}(T),y^g(T)\right>\\=\left<\lambda_k
x^{f_k^1}(T),y^g(T)\right>=\left<\lambda_kC^Tf_k^1,g\right>, \  k
\in \N.
\end{eqnarray*}
Similarly, making use of (\ref{i_2}) for $k=1,\ldots\infty$, $2\leqslant
i\leqslant L_k,$ we obtain
\begin{eqnarray*}
\left<C^T\dot f_k^i,g\right>=\left<x^{\dot
f_k^i}(T),y^g(T)\right>=\left<\dot
x^{f_k^i}(T),y^g(T)\right>\\=\left<\lambda_k
x^{f_k^i}(T)+x^{f_k^{i-1}}(T),y^g(T)\right>=\left<\lambda_kC^Tf_k^i+C^Tf_k^{i-1},g\right>.
\end{eqnarray*}
%For $k=N+1,\ldots,\infty$ we evaluate:
%\begin{eqnarray}
%\left(C^T\dot f_k,g\right)=\left(x^{\dot
%f_k}(T),y^g(T)\right)=\left(\dot
%x^{f_k}(T),y^g(T)\right)\\=\left(\lambda_k
%x^{f_k}(T),y^g(T)\right)=\left(\lambda_kC^Tf_k,g\right)
%\end{eqnarray}

Using (\ref{ct}), one gets the following integral eigenvalue equations for finding   $\lambda_k$ and
$f_k^i$, $1\leqslant i\leqslant L_k$,  $ k \in \N$:
\begin{eqnarray*}
\int_0^Tr(2T-t-\tau)\dot f_k^1(\tau)-\lambda_k
r(2T-t-\tau)f_k^1(\tau)\,d\tau=0,\\
\int_0^Tr(2T-t-\tau)\dot f_k^i(\tau)-\lambda_k
r(2T-t-\tau)f_k^i(\tau)-r(2T-t-\tau)f_k^{i-1}(\tau)\,d\tau=0.
%\\
%\int_0^Tr(2T-t-\tau)\dot f_k(\tau)-\lambda_k
%r(2T-t-\tau)f_k(\tau)\,d\tau=0,\,\,k=N+1,\ldots\infty.
\end{eqnarray*}
Integrating by parts  we finally have:
\begin{eqnarray*}
\int_0^T\dot r(2T-t-\tau) f_k^1(\tau)-\lambda_k
r(2T-t-\tau)f_k^1(\tau)\,d\tau=0,\\
\int_0^T\dot r(2T-t-\tau) f_k^i(\tau)-\lambda_k
r(2T-t-\tau)f_k^i(\tau)-r(2T-t-\tau)f_k^{i-1}(\tau)\,d\tau=0.
%\\
%\int_0^T\dot r(2T-t-\tau) f_k(\tau)-\lambda_k
%r(2T-t-\tau)f_k(\tau)\,d\tau=0,
\end{eqnarray*}
This leads to the following conclusion: the set $\lambda_k,$
$f_k^i$,  $i=1,\ldots,L_k,$ $k \in \N ,$ are eigenvalues and
root vectors of the following generalized eigenvalue problem in
$L_2(0,T)$:
\begin{equation}
\label{gen_ev_prob} \int_0^T\dot r(2T-t-\tau) f(\tau)-\lambda
r(2T-t-\tau)f(\tau)\,d\tau=0.
\end{equation}
Using the same arguments we can deduce that $\overline\lambda_k,$
$g_k^i$, $k=1,\ldots\infty$, $i=1,\ldots,L_k$ are eigenvalues and
root vectors of the eigenvalue problem
\begin{equation}
\label{gen_ev_prob_adj} \int_0^T\overline{\dot r(2T-t-\tau)}
g(\tau)-\lambda \overline{r(2T-t-\tau)}g(\tau)\,d\tau=0.
\end{equation}
We notice that solving (\ref{gen_ev_prob}) and
(\ref{gen_ev_prob_adj}) yields eigenvalues $\lambda_k$,
 their multiplicities $L_k$,
$k \in \N$, and non-normalized functions $f_k^i$ and
$g_k^i$ for which $x^{f_k^i}(T)=\alpha_k^i\phi_{k}^i$,
$y^{g_k^i}(T)=\beta_k^i\psi_{k}^i$, with some (unknown) constants
$\alpha_k^i$, $\beta_k^i$.

Now we describe the algorithm of recovering $a_k^1,\ldots
a_k^{L_k}$, $k \in \N$ (see the representation
(\ref{r_representation})). We normalize the solutions to
(\ref{gen_ev_prob}), (\ref{gen_ev_prob_adj}) by the rule
\begin{equation}
\label{norming} \left<C^T\widetilde f_k^i,\widetilde
g_k^i\right>=1.
\end{equation}
So if $x^{f_k^i}(T)=\phi_k^i$ and $y^{g_k^i}(T)=\psi_k^i$, then
$x^{\widetilde f_k^i}(T)=\alpha_k^i\phi_k^i$ and $y^{\widetilde
g_k^i}(T)=\frac{1}{\alpha_k^i}\psi_k^i$. In the case we define
\begin{eqnarray}
\widetilde b_k^i=\left<y^{\widetilde
g_k^i}(T),b\right>=\int_0^T\overline r(T-\tau)\widetilde
g_k^i(\tau)\,d\tau,\label{b_i}\\
\widetilde d_k^i=\left<x^{\widetilde f_k^i}(T),d\right>=\int_0^T
r(T-\tau)\widetilde f_k^i(\tau)\,d\tau.\label{d_i}
\end{eqnarray}
Then (see (\ref{coeff_repr}))
\begin{equation}
\label{a_k_1} a_k^1=\sum_{i=1}^{L_k}\widetilde b_k^i\widetilde
d_k^i.
\end{equation}
Denote by $\partial$ and $I$ the operator of differentiation and
unitary operator. Bearing in mind (\ref{gen_ev_prob}), which we
rewrite as $C^T\left(\partial-\lambda_k I\right)
f_k^i=C^Tf_k^{i-1}$, we evaluate
\begin{equation*}
%\label{C_ev1}
\left<C^T\left(\partial-\lambda_kI\right)\widetilde
f_k^i,\widetilde g_k^{i-1}\right>=\alpha_k^i\left<C^T
f_k^{i-1},\widetilde
g_k^{i-1}\right>=\frac{\alpha_k^i}{\alpha_k^{i-1}}.
\end{equation*}

So, normalizing the solutions to (\ref{gen_ev_prob}),
(\ref{gen_ev_prob_adj}) by the rule
\begin{equation*}
\left<C^T\left(\partial-\lambda_kI\right)\widehat f_k^i,\widehat
g_k^{i-1}\right>=1,
\end{equation*}
we can define
\begin{eqnarray}
\widehat b_k^i=\int_0^T\overline r(T-\tau)\widehat
g_k^i(\tau)\,d\tau,\label{b_i_h}\\
\widehat d_k^i=\int_0^T r(T-\tau)\widehat
f_k^i(\tau)\,d\tau.\label{d_i_h}
\end{eqnarray}
and calculate $a_k^2=\sum_{i=2}^{L_k}\widehat b_k^i \widehat
d_k^{i-1}$, cf. (\ref{coeff_repr}).

Notice that since $C^T$ commutes with the differentiation, we have
for $l<i$: $\left[C^T\left(\partial-\lambda_k I\right)\right]^l
f_k^i=C^Tf_k^{i-l}$. Then
\begin{equation*}
%\label{C_evl}
\left<\left[C^T\left(\partial-\lambda_kI\right)\right]^l\widetilde
f_k^i,\widetilde g_k^{i-l}\right>=\alpha_k^i\left<C^T
f_k^{i-l},\widetilde
g_k^{i-l}\right>=\frac{\alpha_k^i}{\alpha_k^{i-l}}.
\end{equation*}
Again, normalizing the solutions to (\ref{gen_ev_prob}),
(\ref{gen_ev_prob_adj}) (for $i>l$) by the rule
\begin{equation}
\label{norming_l}
\left<\left[C^T\left(\partial-\lambda_kI\right)\right]^l\widehat
f_k^i,\widehat g_k^{i-l}\right>=1,
\end{equation}
we define $\widehat b_k^i,$ $\widehat d_k^i$ by (\ref{b_i}),
(\ref{d_i}) and evaluate
\begin{equation}
\label{a_k_l} a_k^l=\sum_{i=l}^{L_k}\widehat b_k^i \widehat
d_k^{i-l}.
\end{equation}

We conclude this section with the algorithm for solving the
spectral estimation problem: suppose that we are given with the
function $r\in L_2(0,2T)$ of the form (\ref{r_representation}) and
the family
$\bigcup_{k=1}^\infty\{e^{\lambda_kt},\ldots,t^{L_k-1}e^{\lambda_kt}\}$
is minimal in $L_2(0,T)$. Then to recover $\lambda_k$, $L_k$ and
coefficients of polynomials, one should follow the

{\bf Algorithm}
\begin{itemize}
\item[a)] solve generalized eigenvalue problems
(\ref{gen_ev_prob}), (\ref{gen_ev_prob_adj}) to find $\lambda_k$,
$L_k$ and non-normalized controls.

\item[b)] Normalize $\widetilde f_k^i,$ $\widetilde g_k^i$ by
(\ref{norming}), define $\widetilde b_k^i,$ $\widetilde d_k^i$ by
(\ref{b_i}), (\ref{d_i}) to recover $a_k^1$ by (\ref{a_k_1})(see
(\ref{r_representation}), (\ref{coeff_repr}))

\item[c)] Normalize $\widehat f_k^i,$ $\widehat g_k^i$ by
(\ref{norming_l}), define $\widetilde b_k^i,$ $\widehat d_k^i$ by
(\ref{b_i_h}), (\ref{d_i_h}) to recover $a_k^l$ by
(\ref{a_k_l})(see (\ref{r_representation}), (\ref{coeff_repr}))

\end{itemize}

\section{Continuation of the inverse data for the first order hyperbolic system}

We consider the  initial boundary value problem
\begin{eqnarray}
\frac{\p}{\p t}\begin{pmatrix} u\\
v\end{pmatrix}-\frac{\p}{\p x}\begin{pmatrix} 0&1\\
1&0\end{pmatrix}\begin{pmatrix} u\\
v\end{pmatrix}-\begin{pmatrix} p_{11}& p_{12}\\
p_{21}&p_{22}\end{pmatrix}\begin{pmatrix} u\\
v\end{pmatrix}=0, \quad 0\leqslant x\leqslant 1,\,\, t>0, \label{eq1}\\
u(0,t)=u(1,t)=0,\quad t>0,\label{eq2}\\
\begin{pmatrix} u(x,0)\\
v(x,0)\end{pmatrix}=\begin{pmatrix} d_1(x)\\
d_2(x)\end{pmatrix},\quad 0\leqslant x\leqslant 1\label{eq3}.
\end{eqnarray}
Here   $p_{ij}\in C^1([0,1];\mathbb{C})$ and
  $d_1, d_2 \in L_2(0,1;\mathbb{C})$. We fix some $T>0$
and define
$R(t):=\{v(0,t),v(1,t)\}$, $0\leqslant t\leqslant T.$ The problem
of the recovering unknown potential $p_{ij}$ and initial state
$c_{1,2}$ has been considered in \cite{TY1,TY5}, where the authors
established the uniqueness result for large enough $T.$ The inverse problem by one
measurement for the one-dimensional Schr\"odinger equation has
been considered in \cite{AMR}, and the procedure of the recovering
the potential and the initial state has been proposed. Here we
focus on the problem of the continuation of the inverse data: we
assume that $R(t)$ is known on the interval $(0,T)$ and  recover it
on the whole real axis.

We introduce the notations $B=\begin{pmatrix} 0&1\\
1&0\end{pmatrix}$, $P=\begin{pmatrix} p_{11}& p_{12}\\
p_{21}&p_{22}\end{pmatrix}$, $D=\begin{pmatrix} d_1\\
d_2\end{pmatrix}$  and the operator $A$ acting by the rule
\begin{equation*}
A\varphi=\left(B\frac{d}{dx}+P\right)\varphi,\quad 0\leqslant
x\leqslant 1
\end{equation*}
with the domain
$$
D(A)=\left\{\varphi=\begin{pmatrix} \varphi_1\\
\varphi_2\end{pmatrix}\in H_1(0,1;\mathbb{C}^2)\,|\,
\varphi_1(0)=\varphi_1(1)=0\right\}
$$
%Then we rewrite (\ref{eq1})--(\ref{eq3}) as
%\begin{eqnarray}
%U_t-AU=0, \quad t>0,\label{wave_eq_op}\\
%u(0,t)=u(1,t)=0,\quad  t>0,\label{bc_op}\\
%U(x,0)=C.\label{ic_op}
%\end{eqnarray}
%And the response is given by $R(t)=\{U^2(0,t),U^2(1,t)\}$, $t>0$.

The adjoint operator
\begin{eqnarray*}
A^*\psi=\left(-B\frac{d}{dx}+P^T\right)\psi,\quad 0\leqslant
x\leqslant 1,
\end{eqnarray*}
has the domain
\begin{eqnarray*}
D(A^*)=\left\{\psi=\begin{pmatrix} \psi_1\\
\psi_2\end{pmatrix}\in  H^1(0,1;\mathbb{C}^2)\,|\, \psi_1(0)=\psi_1(1)=0\right\}
\end{eqnarray*}
The spectrum of the operator $A$ has the following structure (see
\cite{TY1,TY5}): $\sigma(A)=\Sigma_1\cup\Sigma_2,$ where
$\Sigma_1\cap\Sigma_2=\emptyset$ and there exists $N_1 \in \N$ such that
\begin{itemize}
\item[1)]$\Sigma_1$ consists of $2N_1-1$ eigenvalues including
algebraical multiplicities

\item[2)] $\Sigma_2$ consists of infinite number of eigenvalues of multiplicity one

\item[3)] Root vectors of $A$ form  a Riesz basis in $ L_2(0,1;\mathbb{C}^2)$.
\end{itemize}

Let  $m$  denote the algebraical multiplicity of eigenvalue
$\lambda,$ and we introduce the notations:
\begin{eqnarray*}
\Sigma_1=\left\{\lambda^i\in\sigma(A), \, m_i\geqslant 2,\,
1\leqslant i\leqslant N\right\},\\
\Sigma_2=\left\{\lambda_n\in\sigma(A), \, \lambda_n \text{ is
simple },\, n\in \mathbb{Z}\right\}.
\end{eqnarray*}
Let $e_1:=\left(1\atop 0\right)$. The root vectors are introduced
by the following way:
\begin{eqnarray*}
\left(A-\lambda^i\right)\phi^i_1 =0,\quad
\left(A-\lambda^i\right)\phi^i_j =\phi^i_{j-1},\quad 2\leqslant
j\leqslant m_i,\\
\phi^i_j(0)=e_1, \,\, \phi^i_j\in D(A),\,\, 1\leqslant j\leqslant
m_i.
\end{eqnarray*}
For the adjoint operator the following equalities are valid:
\begin{eqnarray*}
\left(A^*-\overline\lambda^i\right)\psi^{i}_{m_i} =0,\quad
\left(A^*-\overline\lambda^i\right)\psi^{i}_j
=\psi^{i}_{j+1},\quad 1\leqslant
j\leqslant m_i-1,\\
\psi^{i}_j(0)=e_1, \,\, \psi^{i}_j\in D(A^*),\,\, 1\leqslant
j\leqslant m_i.
\end{eqnarray*}
For the simple eigenvalues we have:
\begin{eqnarray*}
\left(A-\lambda_n\right)\phi_n =0,\quad
\left(A^*-\overline\lambda_n\right)\psi_{n} =0,\\
\phi_n(0)=\psi_{n}(0)=e_1, \,\, \phi_{n}\in D(A),\,\,\psi_{n}\in
D(A^*).
\end{eqnarray*}
Moreover, the following biorthigonality conditions hold:
\begin{eqnarray*}
\left(\phi^i_j,\psi_{n}\right)=0,\quad
\left(\phi_n,\psi^{i}_j\right)=0,\quad
\left(\phi_k,\psi_{n}\right)=0,\\
\left(\phi^i_j,\psi^{k}_l\right)=0,\quad \text{if $i\not= k$ or
$j\not= l$}.
\end{eqnarray*}
Then we set
\begin{eqnarray*}
\rho^i_j=\left(\phi^i_j,\psi^{i}_j\right),\quad i=i\ldots,N,\quad
j=1,\ldots, m_i,\\
\rho_n=\left(\phi_n,\psi_{n}\right),\quad n\in \mathbb{Z},
\end{eqnarray*}
and introduce the spectral data:
\begin{equation*}
S(P)=\left\{\lambda^i,m_i,\rho^i_j\right\}_{1\leqslant i\leqslant
N}^{1\leqslant j\leqslant
m_i}\cup\left\{\lambda_n,\rho_n\right\}_{n\in \mathbb{Z}}
\end{equation*}

We represent the initial state as the series:
\begin{equation}
\label{D_repr} D=\sum_{i=1}^N
\sum_{j=1}^{m_i}d^i_j\phi^i_j(x)+\sum_{n\in
\mathbb{Z}}d_n\phi_n(x).
\end{equation}
We are looking for the solution to (\ref{eq1})--(\ref{eq3}) in the
form
\begin{equation*}
\begin{pmatrix}u\\ v\end{pmatrix}(x,t)=\sum_{i=1}^N
\sum_{j=1}^{m_i}c_j^i(t)\phi_j^i(x)+\sum_{n\in
\mathbb{Z}}c_n(t)\phi_n(x).
\end{equation*}
Using the method of moments we can derive the system of ODe's for
$c^i_j,$ $i\in\{1,\ldots,N\}$, $j\in \{1,\ldots,m_i\}$; $c_n$,
$n\in \mathbb{Z}$ solving which we obtain
\begin{eqnarray*}
c^i_j(t)=e^{\lambda_i t}\left[d^i_j+d^i_{j+1}t+d^i_{j+2}\frac{t^2}{2}+\ldots+d^i_{m_i}\frac{t^{m_i-j}}{(m_i-j)!}\right],\\
c_n(t)=d_ne^{\lambda_n t}.
\end{eqnarray*}
Notice that the response $\{v(0,t),v(1,t)\}$ has a form depicted
in (\ref{unknown_f}):
\begin{eqnarray}
v(0,t)=\sum_{i=1}^N e^{\lambda_i t} a^0_i(t)+\sum_{n\in \mathbb{Z}}e^{\lambda_n t}d_n\left(\phi_n(0)\right)_2,\label{r1}\\
v(1,t)=\sum_{i=1}^N e^{\lambda_i t} a^1_i(t)+\sum_{n\in
\mathbb{Z}}e^{\lambda_n t}d_n\left(\phi_n(1)\right)_2,\label{r2}
\end{eqnarray}
where the coefficients of
$a^0_i(t)=\sum_{k=0}^{m_i-1}\alpha_k^it^k$ are given by
\begin{eqnarray*}
\alpha^i_0=\sum_{l=1}^{m_i}d^i_l\left(\phi^i_{l}(0)\right)_2,\quad \alpha^i_1=\sum_{l=2}^{m_i}d^i_l\left(\phi^i_{l-1}(0)\right)_2,\quad \alpha^i_2=\frac{1}{2}\sum_{l=3}^{m_i}d^i_l\left(\phi^i_{l-2}(0)\right)_2,\\
\ldots,
\alpha^i_k=\frac{1}{(k-1)!}\sum_{l=k+1}^{m_i}d^i_l\left(\phi^i_{l-k}(0)\right)_2,\ldots\quad
\alpha^i_{m_i-1}=\frac{1}{(m_i-1)!}d^i_{m_i}\left(\phi^i_{l}(0)\right)_2.
\end{eqnarray*}
The coefficients $a^1_i(t)$, $i=1,\ldots,N$ are defined by the
similar formulaes.

We introduce the following
\begin{definition}
The state $D\in L_2\left((0,1);\mathbb{C}^2\right)$ is generic if
all the Fourier coefficients in the expansion (\ref{D_repr}) are
not equal to zero.
\end{definition}
We assume below that the initial state $D$ is generic. The meaning
of this restriction is clear -- if the initial state is not
generic, say $d_k=0$ for some $k\in \mathbb{Z}$, the response
$(\ref{r1})$, $(\ref{r2})$ does not contain any information on
$\lambda_k$.

We introduce the notation $U:= \left(u\atop v\right)$ and
consider the dynamical system with the boundary control $f\in
L_2(\mathbb{R}_+)$
\begin{eqnarray*}
U_t-AU=0, \quad 0\leqslant x\leqslant 1,\,\, t>0,\\
u(0,t)=f(t), u(1,t)=0,\quad  t>0, \\
U(x,0)=0.
\end{eqnarray*}
It is not difficult to show that this system is exactly
controllable in time $T \geq 2$. This implies (see \cite{AI}) that
the family
$\bigcup_{i=1}^N\{e^{\lambda_it},\ldots,t^{m_i-1}e^{\lambda_it}\}\cup\{e^{i\lambda_nt}\}_{n\in
\mathbb{Z}}$ forms a Riesz basis in a closure of its linear span
in $L_2((0,T);\mathbb{C})$. Because of this and the fact that each
component of the response $\{v(0,t),v(1,t)\}$ has the form of
(\ref{unknown_f}), we can apply the method from the previous
section and recover $\lambda^i$, $m_i$, coefficients of
polynomials $a^{0,1}_i(t)$ $i=1,\ldots,N,$ $\lambda_n$, $n\in
\mathbb{Z}$. The latter allows one to extend the inverse data
$R(t)$ to all values of $t\in \mathbb{R}$ by formulas (\ref{r1}),
(\ref{r2}). This is important in the solution of the
identification problem, see \cite{TY5}.

%This allows us to make use of the uniqueness results due to
%Trooshin and Yamamoto \cite{TY}. We formulate the final result as
%a theorem.
%\begin{theorem}
%Let us consider the initial boundary value problem
%(\ref{wave_eq_op})--(\ref{ic_op}) with the initial state $D\in
%\left(L_2(0,1)\right)^2$ be generic. Then the observation $R(t)$,
%$t\in (0,2)$ uniquely determines the coefficients $\sigma$,
%$\delta$ and the initial state $D$.
%\end{theorem}

\end{document}